\documentclass[12pt]{amsart}
\usepackage{amsmath}
\usepackage{amsfonts,amssymb,amsthm}
\newtheorem{theorem}{Theorem}[section]
\newtheorem{lemma}[theorem]{Lemma}
\newtheorem{remark}[theorem]{Remark}

\begin{document}
\title{Another simple proof of a theorem of Chandler Davis}
\author{Igor Rivin}
\address{Mathematics Department, Temple University, 1805 N Broad St,
  Philadelphia, PA 19122}
\address{Mathematics Department, Princeton University, Princeton,
NJ 08544}
\email{rivin@math.temple.edu}
\thanks{The author was partially supported by the Department of
  Mathematical Sciences of the National Science Foundation. He would
  like to thank Omar Hijab for helpful conversations.}
\begin{abstract}
In 1957, Chandler Davis proved that  unitarily invariant convex
 functions on the space of hermitian matrices are precisely those
 which are convex and symmetrically invariant on the set of diagonal
 matrices. We give a simple perturbation theoretic proof of this
 result. (Davis' argument was also very short, though based on
 completely different ideas).
\end{abstract}
\maketitle
 Consider an orthogonally invariant function $f$ defined on
the set of $n \times n$ symmetric matrices. Such a function has to
factor through the spectrum:
\begin{equation}
\label{decomp}
f(M) = g \circ \mathbf{\lambda}(M),
\end{equation}
where $g$ is a symmetric function:
\begin{equation}
g(\lambda_1, \dots, \lambda_n) = g(\lambda_{\sigma(1)}, \dots,
\lambda_{\sigma(n)}),
\end{equation}
for any permutation $\sigma.$

In the sequel we shall further assume that $f$ is a $C^2$ convex
function, and under this assumption we shall show that such functions
are precisely those decomposing as per Eq. (\ref{decomp}), with convex
$g.$ The argument and the statement are identical for unitarily
invariant functions of Hermitian matrices; in that setting the theorem
was proved in \cite{davis}, by a completely different argument (Davis
made no regularity assumption, but this is easily dispensed with (see
Remark \ref{regularity})).

To show this, let
$M = P + t Q,$ and let $\tilde{f}_{P, Q}(t) = f(M).$ It is enough to show
that for any symmetric $P, Q,$
$$\frac{d^2\tilde{f}_{P, Q}}{d t^2}(0) > 0.$$

We compute (dropping the subscript, since from now on $P, Q$ do not vary):

\begin{gather}
\frac{d \tilde{f}}{d t} = \sum_{i=1}^n \frac{\partial g}{\partial
  \lambda_i} \dot{\lambda_i}.\\
\frac{d^2 \tilde{f}}{d t^2} =
\sum_{1\leq i, j\leq n} \frac{\partial^2 g}{\partial \lambda_i
  \partial \lambda_j} \dot{\lambda}_i \dot{\lambda_j}
+
\sum_{i=1}^n \frac{\partial g}{\partial \lambda_i} \ddot{\lambda}_i.
\end{gather}

The first sum is positive, since it equals
$$\dot{\mathbf{\lambda}}^t H(g) \dot{\mathbf{\lambda}},$$ and the
Hessian $H(g)$ is positive definite by assumption.
It now suffices to show that the second sum (which can be written as
$\nabla g \cdot \ddot{\mathbf{\lambda}}$) is non-negative.

By continuity, it is sufficient to prove this for $P$ whose spectrum
is simple. By orthogonal invariance, we can compute in a basis where $P$
is diagonal. According to \cite[page 81]{kato}, in that case
$$
\ddot{\lambda}_i = - \sum_{j\neq i} (\lambda_j - \lambda_i)^{-1}
Q_{ij}^2,
$$
and so
$$\nabla f \cdot \ddot{\mathbf{\lambda}} =
\sum_{j > i} Q_{ij}^2 \frac{\frac{\partial g}{\partial \lambda_j} -
    \frac{\partial g}{\partial \lambda_i}}{\lambda_j - \lambda_i}.
$$

The result then follows from the Lemma below.

\begin{lemma}
Let $f$ be a convex function such that $f(x, y) = f(y, x),$ where $y
\neq x.$
Then
$$\frac{\frac{\partial f}{\partial x} - \frac{\partial f}{\partial
    y}}{x - y} > 0.$$
\end{lemma}

\begin{proof}
The conclusion of the lemma is obviously equivalent to:
$$\left(\frac{\partial f}{\partial x} - \frac{\partial f}{\partial
    y}\right)\left(x - y\right) > 0.$$
To show this, consider
$$h(t) = f((1-t) x + t y, t x + (1-t) y).$$
The function $h(t)$ is convex, and $h(0) = h(1).$ This obviously
implies that $\dot{h}(0) \leq 0.$
Since
$$\dot{h}(0) = \left(\frac{\partial f}{\partial x} - \frac{\partial
  f}{\partial y}\right)(y - x),$$ the conclusion follows.
\end{proof}
\begin{remark}
\label{regularity}
Since smooth symmetric convex functions are dense in the set of all symmetric
functions, we have actually shown the result for all convex symmetric
functions, as follows:

\begin{itemize}
\item
We start with a symmetric locally integrable function $g$ on the set
of diagonal matrices.
\item
The function $g$ can be extended (by orthogonal invariance) to the set
of all symmetric matrices, to obtain a function $f.$
\item
Convolving $f$ with a $C^\infty$ Gaussian approximation $D_n$ to the
delta function, then averaging over the orthogonal orbits we obtain a
convex orthogonally invariant function $f_n,$ such that the $\| f -
f_n\| \leq 1/n,$ where $\|\ \|$ denotes the $\mathop{sup}$-norm. The
restriction of $f_n$ to the diagonal matrices is smooth and convex
(since convolving with a positive kernel and averaging both preserve
convexity) function $g_n.$
\item
We apply the smooth argument above to $g_n$ and $f_n,$ to show that
$f_n$ is a convex function.
\item
Since $f$ is a sup-norm limit of convex functions $f_n,$ the result
follows.
\end{itemize}
\end{remark}
\begin{remark}
The computation above can also be extended to certain infinite
dimensional situations. In particular, it is not hard to see that if
$P$ has a compact resolvent and $Q$ is \emph{bounded} then if the
function $g$ is the the spectral zeta function or the logarithm of the
determinant of the operator. In particular, the regularized $\log
\det$ of the Laplacian on a Riemann surface is convex with respect to
bounded perturbations, as are special values of the spectral zeta
function. This is very much \emph{not} true for arbitrary
perturbations: in particular, when $P = Q = \Delta,$ the convexity
h$\log \det$ depends on the sign of the spectral zeta function at $0,$
which is the same as the sign of the Euler characteristic, so that
$\log \det$ is \emph{not} convex with respect to the ``diagonal''
perturbation of the Laplacian whenever the Riemann surfvace has
negative Euler characteristic.
\end{remark}
\bibliographystyle{alpha}

\end{document}